%% file: Beck_Braun_NonHarmonicColorings_AMS_Final.tex
\documentclass[11pt]{amsart}

\usepackage{amsthm,amssymb,amsfonts,amsmath,enumerate,graphicx,mathptmx,microtype}
\newcommand{\Z}{\mathbb{Z}}
\newcommand{\N}{\mathbb{N}}
\newcommand{\R}{\mathbb{R}}

\DeclareMathOperator*{\Ker}{Ker}
\renewcommand{\H}{{\mathcal H}}
\renewcommand{\P}{{\mathcal P}}
\renewcommand{\phi}{\varphi}

\def\Im{\operatorname{Im}}
\def\th{^{\text{th}}}
\def\st{^{\text{st}}}

\def\s{{\boldsymbol s}}

\def\Cz{\R^V}
\def\Co{\R^E}

\newtheorem{theorem}{Theorem}[section]

\newtheorem{proposition}[theorem]{Proposition}
\newtheorem{lemma}[theorem]{Lemma}

\theoremstyle{remark}
\newtheorem{remark}[theorem]{Remark}
\theoremstyle{definition}
\newtheorem{definition}[theorem]{Definition}

\begin{document}


\title{Nowhere-harmonic colorings of graphs}

\author{Matthias Beck}
\address{Department of Mathematics\\
         San Francisco State University\\
         San Francisco, CA 94132}
\email{beck@math.sfsu.edu}
\urladdr{http://math.sfsu.edu/beck}

\author{Benjamin Braun}
\address{Department of Mathematics\\
         University of Kentucky\\
         Lexington, KY 40506--0027}
\email{benjamin.braun@uky.edu}
\urladdr{http://www.ms.uky.edu/~braun/}

\thanks{This research was partially supported by the NSF through grants DMS-0810105 (Beck) and DMS-0758321 (Braun).  The authors would like to thank Tom Zaslavsky and the anonymous referees for their comments and suggestions.}


\date{8 October 2010} 


\begin{abstract}
Proper vertex colorings of a graph are related to its boundary map $\partial_1$, also called its signed vertex-edge incidence matrix.  The vertex Laplacian of a graph, $L=\partial_1 \partial_1^t$, a natural extension of the boundary map, leads us to introduce nowhere-harmonic colorings and analogues of the chromatic polynomial and Stanley's theorem relating negative evaluations of the chromatic polynomial to acyclic orientations. Further, we discuss several examples demonstrating that nowhere-harmonic colorings are more complicated from an enumerative perspective than proper colorings.
\end{abstract}

\maketitle

\section{Introduction}

A classical line of investigation in graph theory is the study of proper graph colorings.  Two well-known enumerative theorems regarding such colorings are the Birkhoff--Whitney theorem establishing the existence of the chromatic polynomial and Stanley's reciprocity theorem, given as Theorems~\ref{chrompoly} and~\ref{graphrec} below.  Proper vertex colorings of a graph are related to its boundary map $\partial_1$, also called its signed vertex-edge incidence matrix; this relationship was exploited by the first author and T. Zaslavsky in \cite{BeckZasIOP} to produce geometric proofs of Theorems \ref{chrompoly} and \ref{graphrec}.  In this paper, the vertex Laplacian of a graph, $L=\partial_1 \partial_1^t=D-A$, a natural extension of the boundary map, leads us to introduce \emph{nowhere-harmonic colorings}.  We prove Theorems \ref{noharmpoly} and \ref{noharmrec}, intriguing analogues of Theorems \ref{chrompoly} and \ref{graphrec}.  Further, we discuss several examples demonstrating that nowhere-harmonic colorings are more complicated from an enumerative perspective than proper colorings.  

Our paper is organized as follows.  In Section~\ref{colorings} we establish notation and review proper graph colorings, including Theorems~\ref{chrompoly} and~\ref{graphrec}.  In Section~\ref{MainResults} we introduce nowhere-harmonic colorings and state our main results.  In Section~\ref{iopsection} we discuss the theory of inside-out polytopes.  In Section~\ref{proofs} we provide proofs of our theorems.  In Section~\ref{examples} we provide the results of some computer calculations and discuss in detail the nowhere-harmonic enumerators for star graphs.


\section{Proper Graph Colorings}\label{colorings}

Let $V=\{x_1,x_2,\ldots,x_n\}$ be a set ordered by subscripts.  
A \emph{graph} $G$ on $V$ is a pair $(V,E)$ where $E$ is a collection of $2$-element subsets of $V$.  As usual, we call the elements of $V$ \emph{vertices} of $G$ and the elements of $E$ \emph{edges} of $G$.  An \emph{$m$-coloring} of $G$ is a labeling of the vertices of $G$ with values from $\{1,\ldots,m\}$ and a \emph{proper $m$-coloring} of $G$ is a labeling of the vertices of $G$ with labels from the set $\{1,\ldots,m\}$ such that adjacent vertices have distinct labels.  We denote by $\chi_G(m)$ the number of proper $m$-colorings of $G$.  The following theorem is the first prerequisite for our investigations.

\begin{theorem}[Birkhoff \cite{Birk}, Whitney \cite{Whitney}]\label{chrompoly}
For every graph $G$ with $n$ vertices, $\chi_G(m)$ is given by a polynomial in $m$ of degree $n$.
\end{theorem}

We call the polynomial $\chi_G(m)$ the \emph{chromatic polynomial} of $G$.  A priori, the value of this polynomial has a combinatorial interpretation only for positive values of $m$; a beautiful theorem providing a combinatorial interpretation for negative integral values of $m$ is the second prerequisite for our investigations.

We refer to the (edge) orientation of $G$ given by directing edge $\{x_i,x_j\}$ as $\{x_i<x_j\}$ whenever $i<j$ as the \emph{standard orientation} of $G$. 
We say an arbitrary orientation $\varepsilon$ of the edges of $G$ is \emph{acyclic} if it does not contain a directed cycle in $G$.  Note that the standard orientation is acyclic.  We may represent any orientation $\varepsilon$ as a sign function on the edges of $G$, i.e., $\varepsilon(e)=1$ if the orientation of $e$ under $\varepsilon$ agrees with the standard orientation and $\varepsilon(e)=-1$ otherwise.
We say a vertex coloring $c$ of $G$ is \emph{compatible} with $\varepsilon$ if for every directed edge $e=\{x_i<x_j\}\in E$ in the standard orientation, either $\varepsilon(e)$ and $c(x_j)-c(x_i)$ have the same sign or $c(x_j)-c(x_i)=0$.

\begin{theorem}[Stanley \cite{StanleyAcyclic}]\label{graphrec}
For every graph $G$ with $n$ vertices and every integer $m> 0$, \[(-1)^{n}\chi_G(-m)=\sum_c\alpha(c),\] where the sum ranges over all $m$-colorings $c$ of $G$ and $\alpha(c)$ is the number of acyclic orientations of $G$ compatible with $c$.
In particular, $(-1)^{n}\chi_G(-1)$ equals the number of acyclic orientations of $G$.
\end{theorem}

One proof of Theorem \ref{graphrec} is via the theory of \emph{inside-out polytopes} developed by the first author and T. Zaslavsky in \cite{BeckZasIOP}.  The theory of inside-out polytopes is a framework for enumerating the points of $\mathbb{Z}^n$ contained in the relative interior of a convex polytope $\P\subset \mathbb{R}^n$ but not in a given affine hyperplane arrangement $\mathcal{H}$, where both the polytope and hyperplane arrangement are rational.  Here a \emph{hyperplane arrangement} is a finite collection of affine hyperplanes in some Euclidean space.  One may view inside-out theory as a template for producing geometric reciprocity theorems for which a combinatorial interpretation may be sought.  In the case of proper graph colorings, geometry enters the picture in the following way.

Given a graph $G$, we denote by $\Cz$ the real vector space with basis given by the vertices of $G$; we denote by $\Co$ the real vector space with basis given by the edges of $G$.  We may then think of a vertex $m$-coloring of $G$ as an integer point in $\Cz$.  It is immediately clear that a proper coloring of a graph $G$ is an integer point in $\Cz \setminus \bigcup \H(G)$, where $\H(G)$ is the hyperplane arrangement induced by $G$, namely \[  \H(G) := \left\{ x_i = x_j : \, \{x_i,x_j\} \in E \right\} .\]  This setup is the starting point of \cite{BeckZasIOP} for studying the chromatic polynomial and acyclic orientations of $G$ from a geometric point of view; we will formulate this in a more general way as follows.

For every graph $G=(V,E)$, there is a natural linear \emph{boundary map} $\partial_1 : \Co \to \Cz$, where $\partial_1(e) = x_j - x_i$ if $e = \{x_i<x_j\}\in E$.
If we represent $\partial_1$ as a matrix, then its transpose $\partial_1^t : \Cz \to \Co$ is the \emph{coboundary map}.  We will liberally use the same notation for a linear transformation and its matrix representation.  If $f: \R^n \to \R^k$ is a linear transformation, then a hyperplane arrangement $\H$ in $\R^k$ induces a hyperplane arrangement $\H^\sharp(f)$ in $\R^n$, the \emph{pullback} of $\H$, as follows: any hyperplane $h \in \H$ is given by a linear functional $\phi$ on $\R^k$, and thus $\phi \circ f$ is a linear functional which, in turn, defines a hyperplane in $\R^n$; the collection of these hyperplanes is $\H^\sharp(f)$.  (Pullbacks were used in \cite{BeckZasIOP} to study antimagic labellings.)  We denote by $\H_0 = \left\{ x_e = 0:e\in E\right\}$ the arrangement of the coordinate hyperplanes in $\Co$ (or in $\Cz$; which space is being referenced will be clear from context).

\begin{lemma}
The proper colorings of $G$ are precisely the positive integer points in $\Cz \setminus \bigcup \H_0^\sharp(\partial_1^t)$.
\end{lemma}

\begin{proof}
We need to show that $\H_0^\sharp(\partial_1^t) = \H(G)$.
The hyperplane $h: x_e = 0$ in $\H_0$ induces the hyperplane $\partial_1^t (\phi) = 0$ in $\H_0^\sharp(\partial_1^t)$, where $\phi$ is the linear form associated with $h$.  But $\partial_1^t \circ \phi$ is the zero matrix except for the $e\th$ row, which is identical to the $e\th$ row of $\partial_1^t$. Thus the hyperplane given by $\partial_1^t (\phi) = 0$ is $x_j - x_i = 0$, where $e = \{x_i,x_j\}$.
\end{proof}

Thus, proper colorings of a graph $G$ arise naturally as positive integer points in $\Cz$ that avoid the hyperplane arrangement $\H_0^\sharp(\partial_1^t)$.  While the combinatorial and geometric importance of the boundary map for $G$ is readily apparent, the relation between the coordinate arrangement $\H_0$ in $\Co$ and $G$ is less apparent.  The importance of the coordinate arrangement can be seen by noting that there is a natural one-to-one correspondence between regions of $\H_0$ and orientations of the edges of $G$, via the following construction:
\begin{itemize}
\item Mark each region of $\H_0$ with a sign vector $\s \in \{ \pm 1 \}^m$, where each sign indicates whether the corresponding coordinate for points in this region is larger or smaller than $0$;
\item Apply $\s$ to the standard orientation of the edges of $G$, reorienting an edge precisely if the corresponding coordinate of $\s$ is negative.
\end{itemize}
This construction can be reversed and thus gives the desired correspondence.  (This idea is not new and goes back to Greene.)  As demonstrated in \cite{BeckZasIOP}, it is this correspondence that is at the heart of Theorem \ref{graphrec}. 

In addition to proper colorings, nowhere-zero flows and tensions on graphs arise from a geometric perspective by studying the boundary map and coordinate hyperplane arrangements.  Enumerative properties of these structures have been the subject of previous papers, such as \cite{BeckZasFlows,breuersanyal,chenwang,kocholflow,kocholtension}.


\section{Nowhere-Harmonic Colorings and Main Results}\label{MainResults}

When in an enumerative situation where inside-out theory applies, one automatically obtains a geometric reciprocity statement that parallels Theorem \ref{graphrec}. However, inside-out theory does not immediately yield a combinatorial interpretation for this reciprocity statement (e.g.,  \cite[Problem 7.2]{BeckZasIOP}).  One of the contributions of this work is to demonstrate that in the case of nowhere-harmonic colorings of graphs, we can find such a combinatorial interpretation via Theorem \ref{noharmrec}.  We will begin by providing an arithmetic definition of nowhere-harmonic colorings.

\begin{definition}
An $m$-coloring $c$ of $G$ is \emph{nowhere harmonic} if for every vertex $v$ of $G$, the value of $c(v)$ is not equal to the average value of $c(w)$ where the average is taken over all $w$ that are neighbors of $v$.  We denote by $\hbar_G(m)$ the number of nowhere-harmonic $m$-colorings of $G$.  We denote the generating function for nowhere-harmonic colorings as \[\bar{H}_G(z):=\sum_{ m \ge 1 } \hbar_G(m)\, z^m.\]
\end{definition}

Two basic properties of nowhere-harmonic colorings and their enumerators are the following.

\begin{proposition}\label{2coloring}
Every graph $G$ admits a nowhere-harmonic $2$-coloring.
\end{proposition}

\begin{proposition}\label{involution}
If $G$ has at least one edge, then $\bar{H}_G(z)$ is divisible by $2z^2$.
\end{proposition}

Given Theorem \ref{chrompoly}, it is natural to ask if the function $\hbar_G(m)$ has nice structure.  One answer is the following; recall first that a \emph{quasipolynomial} is a function $f(m)$ defined on $\Z$ such that there exist $k \in \N$ and polynomials $p_1(m), p_2(m), \dots, p_k(m)$ such that $f(m) = p_j(m)$ whenever $m \equiv j \bmod k$. The minimal such $k$ is the \emph{period} of $f(m)$ and the polynomials $p_1(m), p_2(m), \dots, p_k(m)$ are its \emph{constituents}.

\begin{theorem}\label{noharmpoly}
For every connected graph $G$ with $n$ vertices, there exists a quasipolynomial $f(m)$ of degree $n$ such that $\hbar_G(m)=f(m)$.  Further, $\hbar_G(m)$ has (constant) leading coefficient $1$ and $|\hbar_G(-1)|=2^n-2$.
\end{theorem}

Due to the fact that $\hbar_G(m)$ is a quasipolynomial, $\bar{H}_G(z)$ is a rational function of the form \[\frac{g(z)}{(1-z^p)^{n+1}},\] where $p$ is the period of $\hbar_G(m)$; for justification see \cite[Chapter 4]{StanleyVol1} or \cite[Lemma 3.24]{BeckRobinsCCD}.  As our examples in Section~\ref{examples} indicate, there are subtle and interesting cancellations that occur in these rational forms. 

To state our analogue of Theorem \ref{graphrec}, we need a few preliminary definitions.  Given a graph $G$, a coloring $c$ of $G$, and a vertex $v$ of $G$ of degree $d_v$, we say $c$ is \emph{subharmonic} at $v$ if \[c(v)\leq \frac{1}{d_v}\sum_{\{v,w\}\in E(G)}c(w)\] and we say $c$ is \emph{superharmonic} at $v$ if \[c(v)\geq \frac{1}{d_v}\sum_{\{v,w\}\in E(G)}c(w).\]  We call a $\pm 1$-labeling of the vertices of $G$ a \emph{vertex orientation} of $G$ and say that the vertex orientation is \emph{nonconstant} if it is not the all-ones labeling or the all-negative-ones labeling.  Note that since $G$ has $n$ vertices and two choices of sign for each vertex, there are $2^n-2$ nonconstant vertex orientations of $G$.  Let $\varepsilon$ be a vertex orientation of $G$.  Given a coloring $c$ of the vertices of $G$, we say that $c$ is \emph{compatible} with $\varepsilon$ if $c$ is subharmonic at $v$ when $\varepsilon (v)=1$ and $c$ is superharmonic at $v$ when $\varepsilon (v)=-1$.

\begin{theorem}\label{noharmrec}
For every connected graph $G$ with $n$ vertices and every integer $m> 0$, \[(-1)^{n}\hbar_G(-m)=\sum_c\beta(c),\] where the sum ranges over all $m$-colorings $c$ of $G$ and $\beta(c)$ is the number of nonconstant vertex orientations of $G$ compatible with $c$.
\end{theorem}

As with Theorem \ref{graphrec}, Theorem~\ref{noharmrec} is a consequence of the theory of inside-out polytopes where the role of the boundary map of a graph is played by the vertex Laplacian of a graph.

\begin{definition} Let $G=(V,E)$ be a simple graph with $n$ vertices.  Let the elements of $V$ be ordered as $V=\{x_1,x_2,\ldots,x_n\}$.  The \emph{vertex Laplacian} (also called the \emph{Kirchhoff matrix}) of $G$ is \[L:=D-A,\] where $D$ is the $n\times n$ matrix indexed by the elements of $V$ such that the matrix entry $d_{i,j}$ is the degree of $x_i$ if $i=j$ and $0$ otherwise and $A$ is the $0/1$-valued adjacency matrix of $G$ indexed in the same way as $D$.
\end{definition}

The vertex Laplacian $L$ is connected to the boundary map $\partial_1$ by the well-known equation
\begin{equation}\label{laplacianeq}
  L=\partial_1 \partial_1^t
\end{equation}
(see, e.g., \cite[Chapter 6]{BiggsAGT}),
hence $L$ is a positive semi-definite real symmetric operator and the dimension of the nullspace of $L$ is equal to the number of connected components of $G$.
The identity \eqref{laplacianeq} indicates that the vertex Laplacian is a natural next object to consider from an enumerative point of view, hence our study of nowhere-harmonic colorings.

While the following is well known to algebraic graph theorists, it is helpful for others to compare the vertex Laplacian for a graph $G$ with the usual Laplacian for domains in the complex plane, $\Delta:=\frac{\partial^2}{\partial x^2} + \frac{\partial^2}{\partial y^2}.$  This is an appropriate analogy, as we can write \[\Delta=\mathrm{div} \cdot \mathrm{grad}=\left[ \begin{array}{cc} \frac{\partial}{\partial x} & \frac{\partial}{\partial y}\end{array}\right]\left[ \begin{array}{c} \frac{\partial}{\partial x} \\ \frac{\partial}{\partial y}\end{array}\right],\] which corresponds to the factorization \eqref{laplacianeq}.  Thus, the transpose of the boundary map for a graph is a discrete analogue of the gradient operation for functions on domains while the boundary map for a graph is a discrete analogue of the divergence of a vector field.

A smooth function $f$ on a domain in the plane that satisfies $\Delta f=0$ is called a \emph{harmonic function}.  It should be noted that while harmonic functions on finite graphs are simply described, harmonic functions on infinite graphs have a rich structure and are the subject of active research, see for example \cite{Soardi,LovaszDAFunctions}.  Given a graph $G$ with vertex Laplacian $L$, it is not hard to show that if $c$ is a harmonic function on $G$, i.e., if $Lc=0$, then $c$ is constant on the connected components of $G$ and that if $G$ is connected, then the nullspace of $L$ is spanned by the all-ones vector.  More generally, one has the following combinatorial interpretation of the equation $Lc=b$: if we consider $c$ as a labeling of the vertices of $G$, then the entries in $b$ are the differences between $c(x_i)$ and the average of the values of $c(x_j)$ for all $j$ such that $x_j$ is adjacent to $x_i$.  Thus, we have the following lemma.

\begin{lemma}An $m$-coloring of a graph $G$ is nowhere-harmonic if and only if it corresponds to a positive integer point in $\Cz \setminus \bigcup \H_0^\sharp(L)=\Cz \setminus \bigcup \H_0^\sharp(\partial_1 \partial_1^t)$.
\end{lemma}

The importance of the coordinate arrangement to nowhere-harmonic colorings can be seen by noting that there is a natural one-to-one correspondence between regions of $\H_0$ in $\Cz$ and $\pm 1$-labelings of the vertices of $G$, via the following construction:
\begin{itemize}
\item Mark each region with a sign vector $\s \in \{ \pm 1 \}^m$, where each sign indicates whether the corresponding coordinate for points in this region is larger or smaller than $0$;
\item The resulting vectors $\s$ obviously correspond bijectively to $\pm 1$-labelings of the vertices of~$G$.
\end{itemize}
By specifying a $1$ or $-1$ at each vertex, we are specifying whether that vertex should be a host for sub- or superharmonic behavior for colorings, much as the orientation of an edge restricts the labeling a coloring may have in Theorem \ref{graphrec}.

\begin{remark}
A nowhere-harmonic coloring is not necessarily proper.  It might be interesting to investigate enumerative properties of colorings that are both proper and nowhere harmonic.
\end{remark}


\section{Inside-Out Polytopes}\label{iopsection}

Our proofs of our theorems will involve inside-out theory.  Given a convex polytope $\P\subset \mathbb{R}^n$ and a hyperplane arrangement $\mathcal{H}$, where both the polytope and hyperplane arrangement are rational, i.e. given by linear equations and inequalities with rational coefficients, we call $(\P,\mathcal{H})$ an \emph{inside-out polytope}.  A \emph{region} of $\mathcal{H}$ is a connected component of $\mathbb{R}^n\setminus \bigcup \mathcal{H}$.  A \emph{closed} region is the closure of a region. A \emph{region} of $(\P, \mathcal{H})$ is the nonempty intersection of a region of $\mathcal{H}$ with $\P$.  A \emph{vertex} of $(\P,\mathcal{H})$ is a vertex of any such region. A closed region of $(\P,\mathcal{H})$ is the closure of an open region of $(\P,\mathcal{H})$ and therefore meets the relative interior $\P^\circ$.  The \emph{denominator} of $(\P,\mathcal{H})$ is the minimum $p\in \mathbb{Z}_{\geq 1}$ such that the vertices of the $p$-th dilate of $(\P,\mathcal{H})$ all have integer coordinates.

Given a point $v\in \mathbb{R}^n$, the \emph{multiplicity} of $v$ with respect to $(\P,\mathcal{H})$, denoted $m_{\P,\mathcal{H}}(v)$, is the number of closed regions of $(\P,\mathcal{H})$ containing $v$ if $v\in \P$ and $0$ otherwise.  The \emph{closed Ehrhart quasipolynomial} for $(\P,\mathcal{H})$ is the function defined for $t\in \mathbb{Z}_{\geq 1}$ given by
\[E_{\P,\mathcal{H}}(t):= \sum_{v\in (t^{-1}\mathbb{Z}^n)\cap \P}m_{\P,\mathcal{H}}(v).\]
The \emph{open Ehrhart quasipolynomial} for $(\P,\mathcal{H})$ is the function defined for $t\in \mathbb{Z}_{\geq 1}$ given by
\[E^\circ_{\P^\circ,\mathcal{H}}(t):= \#\left(t^{-1}\mathbb{Z}^n \cap \left[ \P^\circ \setminus \bigcup \mathcal{H} \right]  \right).\]

The main result of inside-out theory is the following.

\begin{theorem}\label{IOThm}{\rm \cite[Theorem 4.1]{BeckZasIOP}} If $(\P,\mathcal{H})$ is a closed, full-dimensional, rational inside-out polytope in $\mathbb{R}^n$, then $E_{\P,\mathcal{H}}(t)$ and $E^\circ_{\P^\circ,\mathcal{H}}(t)$ are quasipolynomials in $t$ that satisfy the reciprocity law $E^\circ_{\P^\circ,\mathcal{H}}(t)=(-1)^nE_{\P,\mathcal{H}}(-t)$, with period equal to a divisor of the denominator of $(\P,\mathcal{H})$, with leading term $\mathrm{vol}(\P)t^n$, and with constant term $E_{\P,\mathcal{H}}(0)$ equal to the number of regions of $(\P,\mathcal{H})$.
\end{theorem}

The main idea behind the proof of Theorem~\ref{IOThm} is that the arrangement $\mathcal{H}$ dissects $\P^\circ$ into disjoint, open, rational polytopes.  Ehrhart theory, see \cite{BeckRobinsCCD} or \cite[Chapter 4]{StanleyVol1}, then implies that the function counting lattice points in dilates of a given region $R$ is a quasipolynomial, called the \emph{Ehrhart quasipolynomial of $R$}, whose degree is equal to the dimension of the region and whose period divides the denomenator of $R$.  By summing up the Ehrhart quasipolynomials of each open region, we recover the counting function $E^\circ_{\P^\circ,\mathcal{H}}(t)$.  One may also compute the function $E^\circ_{\P^\circ,\mathcal{H}}(t)$ via M\"{o}bius inversion on the intersection lattice for the hyperplane arrangement.  The reciprocity law in Theorem~\ref{IOThm} is then a consequence of Ehrhart--Macdonald reciprocity for rational polytopes.


\section{Proofs}\label{proofs}

\begin{proof}[Proof of Proposition \ref{2coloring}]
We may assume that $|V| \ge 2$.
With two colors, the only way that a coloring could avoid being nowhere harmonic would be for some vertex to have the same label as all of its neighbors.  Also, note that we only need show that we can produce such a coloring for connected graphs, as disconnected graphs can be dealt with component by component.  Let $T$ be a rooted spanning tree of $G$ with root $t$.  Label $t$ with $1$.  Label all the neighbors of $t$ with $2$.  Label all of their neighbors with a $1$.  Proceed inductively until all vertices have been exhausted.  When this labeling is complete, every vertex $v$ of $G$ has at least one neighbor with a different label from $v$, hence our coloring is nowhere harmonic.
\end{proof}

\begin{proof}[Proof of Proposition \ref{involution}]
Observe that for any graph $G$ with $n$ vertices and Laplacian $L$, we have $\hbar_G(m-1)=E^\circ_{(0,1)^n,\H_0^\sharp(L)}(m)$.  Let $\P$ denote the inside-out polytope $((0,1)^n,\H_0^\sharp(L))$.  Let $T$ be the invertible affine transformation given by $x\mapsto-x+E$, where $E$ denotes the vector with all entries equal to $1$.  
Since $-x$ is in the cube $[-1,0]^n$ and $E$ is in the kernel of $L$, $T$ takes $\P$ onto itself.  As $T$ is invertible, the bijective correspondence between the regions follows.  Thus, every nowhere-harmonic $m$-coloring $c$ of $G$ has a natural paired coloring $T(c)$, and hence $\bar{H}_G(z)$ is divisible by $2$.  As $G$ has at least one edge, $G$ does not admit nowhere-harmonic $1$-colorings; thus, the first term in $\bar{H}_G(z)$ is the $z^2$ term and our claim is verified.
\end{proof}

\begin{proof}[Proof of Theorem \ref{noharmpoly}] Given that for any graph $G$ with $n$ vertices and Laplacian $L$, we have $\hbar_G(m-1)=E^\circ_{(0,1)^n,\H_0^\sharp(L)}(m)$, Theorem~\ref{IOThm} immediately implies the quasipolynomiality and degree results for $\hbar_G(m)$.  To show that $|\hbar_G(-1)|=2^n-2$, we need to prove that the number of regions of $([0,1]^n,\mathcal{H}_0^\sharp(L))$ is $2^n-2$.  This is a straightforward exercise in linear algebra, as follows.

As $L$ is a positive semi-definite real symmetric operator, the spectral theorem implies that $L$ has only nonnegative real eigenvalues and an orthogonal basis of eigenvectors.  Since $G$ is connected, we can conclude that the all-ones vector spans the kernel of $L$, and that $L$ is invertible on $\Ker (L)^\perp$.  Thus, each hyperplane $K$ 
in the arrangement $\mathcal{H}_0^\sharp (L)$ is spanned by $\Ker (L)$ and $K \cap \Ker (L)^\perp$.  The intersection of $K$ with $\Ker (L)^\perp$ is the inverse image of one of the coordinate hyperplanes in $\Cz$ restricted to $\Im (L)$.  As $L$ restricted to $\Ker (L)^\perp$ is invertible with positive eigenvalues and orthogonal eigenvectors, we can conclude that the regions in $([0,1]^n,\mathcal{H}_0^\sharp(L))$ are in one-to-one correspondence with the regions in $\Im (L)\cap \mathcal{H}_0$.  There are exactly $2^n-2$ such regions, since $\Im (L)$ intersects all the open regions defined by the coordinate arrangement except the positive and negative orthants.
\end{proof}

\begin{proof}[Proof of Theorem \ref{noharmrec}]
Our goal is to interpret the geometric reciprocity portion of Theorem \ref{IOThm} combinatorially in this setting.  From the proof of Theorem \ref{noharmpoly}, we know that \[E_{[0,1]^n,\mathcal{H}_0^\sharp (L)}(m)=(-1)^n \hbar_G(-(m+1))\] and that the regions of $([0,1]^n,\mathcal{H}_0^\sharp (L))$ correspond bijectively with nonconstant vertex orientations of $G$.  To finish the proof, note that $E_{[0,1]^n,\mathcal{H}_0^\sharp (L)}(m)$, by definition, sums over the colorings described in the theorem and the multiplicity of each coloring $c$ in the sum is exactly $\beta(c)$.
\end{proof}


\section{Examples}\label{examples}

Using the {\tt IOP} package \cite{vanherick}, we computed sample generating functions for nowhere-harmonic quasipolynomials among four classes of graphs: path graphs, cycle graphs, complete graphs, and star graphs.  In addition to our computations, we provide in this section a detailed analysis of nowhere-harmonic colorings of star graphs to show how polyhedral geometry can illuminate the structure of these generating functions.

In the examples that follow, though $\bar{H}_G(z)$ is a rational function of the form $\frac{g(z)}{(1-z^p)^{n+1}}$, where $p$ is the period of $\hbar_G(m)$, this is not necessarily the reduced form of the quotient.  We have chosen to represent these rational functions in both unreduced and reduced form to emphasize both $p$ and its divisors.  This will prove particularly useful when we analyze nowhere harmonic functions of star graphs.  The numerators of these expressions remain quite mysterious.

\begin{remark}
Recall from the paragraph preceding Theorem~\ref{noharmpoly} that the period of a quasipolynomial $f(m)$ is the number of constituent polynomials for $f(m)$.  
For a rational polytope $Q$ with vertex coordinates $\frac{r_i}{q_i}$, where the index set $I$ runs through the coordinates of all vertices, one might hope that the period of the Ehrhart quasipolynomial for $Q$ is $\mathrm{lcm}\{q_i:i\in I\}$ (which is the maximal possible period, see, e.g., \cite[Section 3.7]{BeckRobinsCCD}).
If this were the case, then by studying the exponent(s) in the denominator of a rational form of $\sum_{m=0}^\infty f(m)z^m$ one might gain information about the denominators of the vertex coordinates.
However, our hope is ruined by the observation that for general rational polytopes the period of an Ehrhart quasipolynomial can be smaller than $\mathrm{lcm}\{q_i:i\in I\}$; as a result, one must exhibit some caution when drawing numerical conclusions about the vertex denominators from the exponents of the denominators of the algebraic expressions.  
Despite this occasional pathological behavior, when studying (inside-out) polytopes arising from combinatorial constructions it is often the case that the denominators of the relevant vertices are reflected in a reasonable way in the denominators of rational generating functions, as we will see later in this section.
\end{remark}

\subsection{Path Graphs}

The \emph{path graph} on $n$ vertices has the edges $\{j, j+1\}$ for $1 \le j < n$.
We denote its generating function for nowhere-harmonic $m$-colorings by $P_n(z)$.

\begin{align*}
  P_3(z) &= \frac{2z^2 \left( 1+2z+3z^2 \right) }{(1-z^2)(1-z)^3} \\
         &= \frac{ 2 z^2 + 10 z^3 + 24 z^4 + 32 z^5 + 22 z^6 + 6 z^7 }{ (1-z^2)^4 } \\ \\
\end{align*}
\begin{align*}
  P_4(z) &= \frac{ 2 z^2 \left( 2 + 10 z + 31 z^2 + 40 z^3 + 36 z^4 + 18 z^5 + 7 z^6 \right) }{ (1 - z)^2 (1 - z^2)^2 (1 - z^3) } \\ 
         &= \frac{ \left\{
              \begin{aligned}
              4 z^2 + 28 z^3 + 122 z^4 + 340 z^5 + 786 z^6 + 1558 z^7 + 2794 z^8 + 4550 z^{9} + 6794 z^{10} \\
              + 9430 z^{11} + 12210 z^{12} + 14862 z^{13} + 16962 z^{14} + 18246 z^{15} + 18510 z^{16} \\
              + 17718 z^{17} + 15990 z^{18} + 13554 z^{19} + 10814 z^{20} + 8066 z^{21} + 5614 z^{22} \\
              + 3602 z^{23} + 2118 z^{24} + 1130 z^{25} + 530 z^{26} + 214 z^{27} + 64 z^{28} + 14 z^{29} 
              \end{aligned} \right\}
            }{ (1-z^6)^5 } \\ \\
  P_5(z) &= \frac{ 2 z^2 \left\{
              \begin{aligned}
                3 + 30 z + 173 z^2 + 480 z^3 + 979 z^4 + 1456 z^5 + 1740 z^6\\ + 1586 z^7 + 1185 z^8 + 648 z^9 + 273 z^{10} + 72 z^{11} + 15 z^{12} 
              \end{aligned} \right\}
            }{ (1 - z) (1 - z^2)^2 (1 - z^3)^2 (1 - z^4) } \\ 
         &= \frac{ \left\{
              \begin{aligned}
                6 z^2 + 66 z^3 + 424 z^4 + 1516 z^5 + 4310 z^6 + 10098 z^7 + 21092 z^8 + 39816 z^{9} \\+ 70326 z^{10} + 116810 z^{11} + 185604 z^{12} + 283092 z^{13} + 418382 z^{14} + 600078 z^{15} \\+ 839332 z^{16} + 1145352 z^{17} + 1528790 z^{18} + 1997262 z^{19} + 2558320 z^{20} \\+ 3215964 z^{21} + 3971790 z^{22} + 4824022 z^{23} + 5765976 z^{24} + 6787584 z^{25} \\+ 7871332 z^{26} + 8997876 z^{27} + 10139176 z^{28} + 11269592 z^{29} + 12355844 z^{30} \\+ 13370100 z^{31} + 14279040 z^{32} + 15056976 z^{33} + 15675972 z^{34} + 16116804 z^{35} \\+ 16361856 z^{36} + 16403112 z^{37} + 16238436 z^{38} + 15873756 z^{39} + 15323840 z^{40} \\+ 14606000 z^{41} + 13746484 z^{42} + 12770172 z^{43} + 11709592 z^{44} + 10592568 z^{45} \\+ 9452772 z^{46} + 8316652 z^{47} + 7213800 z^{48} + 6164832 z^{49} + 5190070 z^{50} \\+ 4301034 z^{51} + 3506480 z^{52} + 2809692 z^{53} + 2210566 z^{54} + 1705914 z^{55} \\+ 1289372 z^{56} + 953736 z^{57} + 688806 z^{58} + 485522 z^{59} + 332604 z^{60} \\+ 221220 z^{61} + 141614 z^{62} + 87030 z^{63} + 50588 z^{64} + 27688 z^{65} + 13846 z^{66} \\+ 6294 z^{67} + 2424 z^{68} + 780 z^{69} + 174 z^{70} + 30 z^{71}
              \end{aligned} \right\}
            }{ (1 - z^{12})^6 } \\ \\
  P_6(z) &= \frac{ 2 z^2 \left\{
    \begin{aligned}
   & 5 + 82 z + 748 z^2 + 3326 z^3 + 10964 z^4 + 27767 z^5 + 59103 z^6 + 106774 z^7 \\ & + 171025 z^8 + 242771 z^9 + 312401 z^{10} + 363363 z^{11} + 386862 z^{12} \\ & + 374277 z^{13} + 331925 z^{14} + 266383 z^{15} + 194327 z^{16} + 126327 z^{17} \\ & + 73219 z^{18} + 36387 z^{19}  + 15522 z^{20} + 5208 z^{21} + 1372 z^{22} + 231 z^{23} + 31 z^{24}
    \end{aligned}
    \right\}}{ (1 - z^2) (1 - z^3)^2 (1 - z^4)^2 (1 - z^5) (1 - z^6) } 
\end{align*}
In unreduced form, $P_6(z)$ has a denominator of $(1-z^{ 60 })^7$ (and a numerator too large to display on this page).

\subsection{Cycles}

The \emph{cycle} on $n$ vertices has the edges $\{j, j+1\}$, for $1 \le j < n$, and $\{ 1, n \}$.
We denote its generating function for nowhere-harmonic $m$-colorings by $C_n(z)$.

\begin{align*}
  C_3(z) &= \frac{ 6 z^2 (1 + z^2) }{ (1 - z)^3 (1 - z^2) } \\ 
         &= \frac{ 6 z^2 + 18 z^3 + 24 z^4 + 24 z^5 + 18 z^6 + 6 z^7}{ (1 - z^2)^4 } \\ \\
  C_4(z) &= \frac{ 2 z^2 (3 + 11 z + 41 z^2 + 42 z^3 + 29 z^4 + 11 z^5 + 7 z^6) }{ (1 - z)^2 (1 - z^2)^2 (1 - z^3) } \\ 
         &= \frac{ \left\{
    \begin{aligned}
6 z^2 + 34 z^3 + 156 z^4 + 412 z^5 + 930 z^6 + 1798 z^7 + 3170 z^8 + 5078 z^{9} + 7426 z^{10} + 10150 z^{11} \\+ 12930 z^{12} + 15510 z^{13} + 
 17406 z^{14} + 18426 z^{15} + 18426 z^{16} + 17358 z^{17} + 15414 z^{18} \\+ 12834 z^{19} + 10086 z^{20} + 7394 z^{21} + 5046 z^{22} + 3170 z^{23} + 1830 z^{24} + 962 z^{25} + 436 z^{26} \\+ 172 z^{27} + 50 z^{28} + 14 z^{29} 
    \end{aligned}
    \right\} }{ (1 - z^6)^5 } \\ \\
  C_5(z) &= \frac{ 10 z^2 \left\{
    \begin{aligned}
    1 + 7 z + 38 z^2 + 75 z^3 + 148 z^4 + 201 z^5 + 265 z^6 + 256 z^7\\ + 271 z^8 + 195 z^9 + 152 z^{10} + 75 z^{11} + 34 z^{12} + 7 z^{13} + 3 z^{14}   
    \end{aligned}
    \right\} }{ (1 - z)^2 (1 - z^2) (1 - z^3) (1 - z^4) (1 - z^6) } \\ 
         &= \frac{ \left\{
    \begin{aligned}
10 z^2 + 90 z^3 + 560 z^4 + 1860 z^5 + 5110 z^6 + 11650 z^7 + 23940 z^8 + 44460 z^{9} + 77770 z^{10} \\+ 127890 z^{11} + 201780 z^{12} + 305660 z^{13} + 449410 z^{14} + 641190 z^{15} + 892900 z^{16} \\+ 1213080 z^{17} + 1612870 z^{18} + 2098950 z^{19} + 2679000 z^{20} + 3356280 z^{21} + 4131570 z^{22} \\+ 5002710 z^{23} + 5961720 z^{24} + 6998080 z^{25} + 8092300 z^{26} + 9225540 z^{27} \\+ 10367160 z^{28} + 11493000 z^{29} + 12567460 z^{30} + 13564820 z^{31} + 14449920 z^{32} \\+ 15199320 z^{33} + 15784460 z^{34} + 16188180 z^{35} + 16392960 z^{36} + 16393080 z^{37} \\+ 16187660 z^{38} + 15783660 z^{39} + 15198400 z^{40} + 14448720 z^{41} + 13563380 z^{42} \\+ 12566060 z^{43} + 11491560 z^{44} + 10365840 z^{45} + 9224300 z^{46} + 8091180 z^{47} \\+ 6997320 z^{48} + 5960800 z^{49} + 5002250 z^{50} + 4131330 z^{51} + 3356120 z^{52} + 2679060 z^{53} \\+ 2099270 z^{54} + 1613130 z^{55} + 1213500 z^{56} + 893340 z^{57} + 641610 z^{58} + 449850 z^{59} \\+ 306060 z^{60} + 202220 z^{61} + 128210 z^{62} + 78030 z^{63} + 44700 z^{64} + 24120 z^{65} + 11750 z^{66} \\+ 5230 z^{67} + 1920 z^{68} + 600 z^{69} + 130 z^{70} + 30 z^{71}
    \end{aligned}
    \right\} }{ (1 - z^{12})^6 } \\ \\
  C_6(z) &= \frac{ 2 z^2 \left\{
    \begin{aligned}
  &  10 + 108 z + 965 z^2 + 4022 z^3 + 12814 z^4 + 31491 z^5 + 65668 z^6 + 115971 z^7 \\ & + 182775 z^8 + 254824 z^9 + 323136 z^{10} + 369837 z^{11} + 388256 z^{12} + 369420 z^{13} \\ & + 322890 z^{14} + 254413 z^{15} + 182376 z^{16} + 115872 z^{17} + 65683 z^{18} + 31509 z^{19} \\ & + 13018 z^{20} + 4109 z^{21} + 1034 z^{22} + 168 z^{23} + 31 z^{24}
    \end{aligned}
    \right\} }{ (1 - z^2) (1 - z^4)^2 (1 - z^3)^2 (1 - z^5) (1 - z^6) } 
\end{align*}
In unreduced form, $C_6(z)$ has a denominator of $(1-z^{ 60 })^7$.

\subsection{Complete Graphs}

The \emph{complete graph} on $n$ vertices has the edges $\{ i, j \}$ for $1 \le i < j \le n$.
We denote its generating function for nowhere-harmonic $m$-colorings by $K_n(z)$.

\begin{align*}
  K_3(z) &= C_3(z) = \frac{ 6 z^2 (1 + z^2) }{ (1 - z)^3 (1 - z^2) } = \frac{ 6 z^2 + 18 z^3 + 24 z^4 + 24 z^5 + 18 z^6 + 6 z^7}{ (1 - z^2)^4 }  \\ \\
  K_4(z) &= \frac{ 2 z^2 (7 + 12 z + 17 z^2 + 29 z^3 + 7 z^5) }{ (1 - z)^3 (1 - z^2) (1 - z^3) } \\
         &= \frac{ \left\{
              \begin{aligned}
                14 z^2 + 66 z^3 + 204 z^4 + 524 z^5 + 1098 z^6 + 2070 z^7 + 3514 z^8 + 5478 z^{9} + 7914 z^{10} + 10582 z^{11} \\+ 13338 z^{12} + 15750 z^{13} + 17526 z^{14} + 18378 z^{15} + 18162 z^{16} + 17022 z^{17} + 14958 z^{18} + 12402 z^{19} \\+ 9646 z^{20} + 7026 z^{21} + 4782 z^{22} + 2962 z^{23} + 1710 z^{24} + 882 z^{25} + 404 z^{26} + 156 z^{27} + 42 z^{28} + 14 z^{29}
              \end{aligned}
            \right\} }{ (1 - z^6)^5 } \\ \\
  K_5(z) &= \frac{ 10 z^2 (3 + 7 z + 28 z^2 + 19 z^3 + 50 z^4 + 5 z^5 + 32 z^6 - 3 z^7 + 3 z^8) }{ (1 - z)^4 (1 - z^3) (1 - z^4) } \\
         &= \frac{ \left\{
              \begin{aligned}
16336800z^{37}+11207800z^{44}+12295020z^{43}+13310300z^{42}+16081620z^{38} \\
+16388640z^{36}+14223840z^{41}+15006280z^{40}+15632900z^{39}+15882300z^{34} \\
+16234980z^{35}+4470z^{67}+70130z^{63}+10310z^{66}+116370z^{62}+21240z^{65} \\
+705130z^{15}+90390z^{10}+39940z^{64}+6730z^{6}+1660z^{68}+500z^{69}+860z^4 \\
+29060z^8+90z^{70}+30z^{71}+8377740z^{26}+3189100z^{52}+415050z^{59} \\
+595470z^{58}+832100z^{57}+2536320z^{53}+3939590z^{51}+1135540z^{56}+1515330z^{55} \\
+1980050z^{54}+5722800z^{49}+4786230z^{50}+6738720z^{48}+7818300z^{47}+8939940z^{46} \\
+2640z^5+15342520z^{33}+14635400z^{32}+13784340z^{31}+9511900z^{27} \\
+7273680z^{25}+12816580z^{30}+11761200z^{29}+10649720z^{28}+5244870z^{23}+6223200z^{24} \\
+10079480z^{45}+190z^3+146610z^{11}+52900z^9+4351650z^{22} \\
+3552340z^{21}+2850380z^{20}+973940z^{16}+497850z^{14}+2245950z^{19}+1735870z^{18} \\
+281280z^{60}+1314600z^{17}+228000z^{12}+342120z^{13}+14730z^7+184440z^{61}+30z^2
              \end{aligned}
            \right\}  }{ (1 - z^{12})^6 } \\ \\
  K_6(z) &= \frac{ 2 z^2 \left\{
    \begin{aligned}
   & 31 + 179 z + 720 z^2 + 1469 z^3 + 2400 z^4 + 3510 z^5 + 3511 z^6 \\ & + 3841 z^7 + 2520 z^8 + 2220 z^9 + 569 z^{10} + 630 z^{11} - 31 z^{12} + 31 z^{13}
    \end{aligned}
    \right\} }{ (1 - z)^4 (1 - z^3) (1 - z^4) (1 - z^5) } 
\end{align*}
In unreduced form, $K_6(z)$ has a denominator of $(1-z^{ 60 })^7$.

\subsection{Star Graphs}

The \emph{star graph on $n$ vertices} is the complete bipartite graph $K_{1,n-1}$, i.e., the edges are $\left\{\{1,j\}: \, 2 \le j \le n \right\}$.  
We denote its generating function for nowhere-harmonic $m$-colorings by $S_n(z)$.

\begin{align*}
S_3(z) &= P_3(z) = \frac{2z^2(3z^2+2z+1)}{(1-z)^3(1-z^2)} = \frac{ 2 z^2 + 10 z^3 + 24 z^4 + 32 z^5 + 22 z^6 + 6 z^7 }{ (1-z^2)^4 } \\ \\
S_4(z) &= \frac{2z^2(7z^5+12z^4+26z^3+17z^2+9z+1)}{(1-z)^3(1-z^2)(1-z^3)} \\ 
       &= \frac{ \left\{
              \begin{aligned}
2 z^2 + 24 z^3 + 102 z^4 + 308 z^5 + 720 z^6 + 1458 z^7 + 2644 z^8 + 4350 z^9 + 6576 z^{10} \\+ 9178 z^{11} + 11988 z^{12} + 14670 z^{13} + 16860 z^{14} + 18234 z^{15} + 18576 z^{16} + 17886 z^{17} \\+ 16200 z^{18} + 13806 z^{19} + 11044 z^{20} + 8274 z^{21} + 5784 z^{22} + 3718 z^{23} + 2196 z^{24} \\+ 1170 z^{25} + 554 z^{26} + 222 z^{27} + 66 z^{28} + 14 z^{29}
              \end{aligned}
            \right\} }{ (1 - z^6)^5 } \\ \\
S_5(z) &= \frac{2z^2(15z^9+50z^8+177z^7+265z^6+323z^5+309z^4+187z^3+95z^2+18z+1)}{(1-z)^3(1-z^2)(1-z^3)(1-z^4)} \\ 
       &= \frac{ \left\{
              \begin{aligned}
2 z^2 + 42 z^3 + 312 z^4 + 1224 z^5 + 3626 z^6 + 8802 z^7 + 18696 z^8 + 35976 z^9 + 64178 z^{10} \\+ 107778 z^{11} + 172416 z^{12} + 264888 z^{13} + 393350 z^{14} + 567198 z^{15} + 796488 z^{16} \\+ 1091616 z^{17} + 1462238 z^{18} + 1917270 z^{19} + 2463768 z^{20} + 3106608 z^{21} + 3847766 z^{22} \\+ 4685958 z^{23} + 5615520 z^{24} + 6626448 z^{25} + 7703348 z^{26} + 8825556 z^{27} + 9968208 z^{28} \\+ 11102880 z^{29} + 12199748 z^{30} + 13227588 z^{31} + 14156016 z^{32} + 14956128 z^{33} \\+ 15601940 z^{34} + 16071108 z^{35} + 16347168 z^{36} + 16419456 z^{37} + 16285228 z^{38} \\+ 15949740 z^{39} + 15425328 z^{40} + 14730768 z^{41} + 13889116 z^{42} + 12927708 z^{43} \\+ 11875920 z^{44} + 10764528 z^{45} + 9624268 z^{46} + 8485116 z^{47} + 7374048 z^{48} + 6315120 z^{49} \\+ 5327338 z^{50} + 4424418 z^{51} + 3615096 z^{52} + 2903448 z^{53} + 2289970 z^{54} + 1771578 z^{55} \\+ 1342728 z^{56} + 995736 z^{57} + 721498 z^{58} + 509850 z^{59} + 350688 z^{60} + 233928 z^{61} \\+ 150574 z^{62} + 92886 z^{63} + 54408 z^{64} + 29904 z^{65} + 15142 z^{66} + 6894 z^{67} + 2712 z^{68} \\+ 864 z^{69} + 190 z^{70} + 30 z^{71}
              \end{aligned}
            \right\} }{ (1 - z^{12})^6 } \\ \\
S_6(z) &= 
\frac{
2z^2\, \left\{
\begin{aligned}
& 31z^{14}+180z^{13}+934z^{12}+2174z^{11}+4039z^{10}+5980z^9+7256z^8+ \\
& 7417z^7+6386z^6+4615z^5+2629z^4+1189z^3+324z^2+45z+1
\end{aligned} \right\}
 }
 {(1-z)^3(1-z^2)(1-z^3)(1-z^4)(1-z^5)}
\end{align*}
In unreduced form, $S_6(z)$ has a denominator of $(1-z^{ 60 })^7$.

\begin{remark}
Observe that for all the graphs we have computed, the denominators of $\bar{H}_G(z)$ in unreduced form are identical when the graphs have the same number of vertices.  
It would be interesting to determine if this is a general pattern for all graphs or a consequence of our examples being in low dimensions. 
\end{remark}

Note that a star graph on $n$ vertices is a tree with $n-1$ leaves and a single vertex of degree $n-1$.  On the leaves of $K_{1,n-1}$, to satisfy the nonharmonic condition requires only that the labels on the leaves do not equal the label on the center of the star, vertex $1$.  Therefore, in the case of $K_{1,n-1}$, every nowhere-harmonic coloring is actually a proper coloring (but not vice versa). $K_{1,n-1}$ is a simple graph to look at from a nowhere-harmonic point of view, in the sense that after removing the first row of the matrix $L_n$, 
one is left with the matrix corresponding to the transpose of the boundary map for $K_{1,n-1}$, which is totally unimodular.  Thus, by removing the hyperplane corresponding to the first row of $L_n$, we are left with an inside-out polytope whose open Ehrhart polynomial is actually a polynomial, as all the vertices of the regions are integral.  We will see that the addition of this single hyperplane creates regions with large vertex denominators leading to these complicated, though not completely chaotic, generating functions for nowhere-harmonic $m$-colorings.

There are a variety of interesting patterns to observe regarding the functions $S_n(z)$.  As anticipated by Proposition~\ref{involution}, their numerator polynomials are of the form $2z^2f_n(z)$ for some polynomial $f_n(z)$.  Also, the leading terms of the $f_n(z)$'s are of the form $2^{n-1}-1$, which essentially counts regions in the associated inside-out polytope.  A mystery to us is that the coefficients of the second highest terms of the $f_n(z)$'s in reduced form appear to be of the form $3^{n-1}-2\cdot 2^{n-1}+1=2S(n,3)$, where $S(n,3)$ is the Stirling number of the second kind \cite[sequence $R_{10}$]{LaHaye}.  All the denominators of the generating functions in our examples are of the form $(1-z^{n-1})\cdots(1-z^2)(1-z)^3$, though we have not been able to prove that this will always be the case.

Our goal in the remainder of this section is to investigate the inside-out polytopes $\P_n:=((0,1)^n,\mathcal{H}_0^\sharp (L_n))$, where $L_n$ is the vertex Laplacian for $K_{1,n-1}$, with two purposes.  First, through a careful analysis of the regions of $\P_n$, we will be able to produce an alternative expression for the generating function $S_n(z)$.  Second, our analysis will also yield some insight regarding the form of the denominator of $S_n(z)$.

Fix $n\geq 3$.  Observe that the vertex Laplacian for $K_{1,n-1}$ is a block matrix of the form
\[L_n = \bordermatrix{\text{}&1&2&3&\ldots &n\cr
                1&n-1 & -1 & -1 & \cdots & -1\cr
                2&-1 & 1 & 0 & \cdots & 0\cr
                3&-1 & 0 & 1 & \cdots & 0\cr
                \vdots&\vdots & \vdots  & \vdots  &\ddots & \vdots\cr
				n&-1 & 0 & 0 & \cdots & 1 },\]
where the rows and columns of $L_n$ are linearly indexed by $\{1,2,\ldots,n\}$ as shown.  Each region $R_\varepsilon$ of the inside-out polytope $\P_n$ is determined by a nonconstant vertex orientation $\varepsilon$ of $K_{1,n-1}$, with the constraint matrix for $R_\varepsilon$ given by
\[\left( \begin{array}{c}
D(\varepsilon)L_n\\
I_n \\
-I_n \\
\end{array} \right),\]
where $D(\varepsilon)$ denotes the diagonal matrix with $\varepsilon$ on its diagonal; thus, $D(\varepsilon)L_n$ is the matrix obtained by multiplying each row of the Laplacian matrix by a $-1$ or $1$ in accordance with $\varepsilon$.  For example, with $L_3$ and the vertex orientation $\{ 1\mapsto -1,2\mapsto 1,3\mapsto -1\}$, the corresponding region is
\[\left( \begin{array}{rrr}
-2 & 1 & 1 \\
-1 & 1 & 0 \\
1 & 0 & -1 \\
1  &  0& 0 \\
0 & 1 & 0 \\
0 & 0 & 1 \\
-1 & 0 & 0 \\
0 & -1 & 0 \\
0 & 0 & -1 
\end{array} \right)
\left(\begin{array}{c}
x_1\\
x_2\\
x_3
\end{array}\right)
<
\left( \begin{array}{r}
0 \\
0 \\
0 \\
1 \\
1 \\
1 \\
0 \\
0 \\
0 \\
\end{array}\right).\]

\begin{figure}[ht]
\begin{center}
 \input{P3.pstex_t}
\end{center}
\caption{The inside-out polytope $\P_3$.}
\label{P3}
\end{figure}
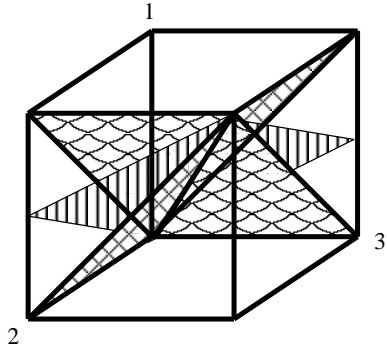

Figure~\ref{P3} contains a picture of $\P_3$, where the axes are labeled $1,2,3$ in accordance with the vertices.  Two features of $\P_3$ are immediately apparent: first, there is a transformation taking the half of $\P_3$ ``above'' the hyperplane $2x_1 + x_2 + x_3 = 0$ to the half of $\P_3$ ``below'' the hyperplane; this is the transformation $T$ from Proposition~\ref{involution}.  Second, there are only two types of regions, one being pyramids over squares.  We can formulate precise analogues of these features for general $\P_n$.

We will from now on consider only regions of $\P_n$ corresponding to vertex orientations of $K_{1,n-1}$ where $1\mapsto -1$.  We will call these regions \emph{negative regions} of $\P_n$.  Our next goal is to produce an alternative formula for $S_n(z)$ by utilizing symmetries of $K_{1,n-1}$.  The automorphism group of $K_{1,n-1}$ consists of $\mathfrak{S}_{n-1}$, the permutations of the $n-1$ leaves of the star; $\mathfrak{S}_{n-1}$ induces a permutation transformation of $\mathbb{R}^n$ under which $\P_n$ has its negative regions permuted.  Two negative regions are in the same orbit of this group action if they correspond to vertex orientations with the same number of $-1$'s.  It follows that there are exactly $n-1$ orbits of this $\mathfrak{S}_{n-1}$-action and, for $1\leq j \leq n-1$, the $j\th$ orbit is represented by the region $\P_n^j$ defined by negating the first $j$ rows of $L_n$, i.e., the region for the vertex orientation $\{1\mapsto -1,2\mapsto -1,\ldots,j\mapsto -1,j+1\mapsto 1,\ldots,n\mapsto 1\}$.  We may now compute $S_{n}(z)$ by summing the open Ehrhart series for the $\P_n^j$'s weighted by the number of elements in each respective orbit, obtaining the following proposition.

\begin{proposition}\label{symmetry}
Let $\P_n^j$ be defined as just described, and let \[S_n^j(z):=\sum_{m=1}^\infty \#\left(m^{-1}\mathbb{Z}^n \cap  \P_n^j \right) z^m.\]  The following equality holds:
\[S_{n}(z)=2\sum_{j=1}^{n-1}\binom{n-1}{j-1}S_n^{j}(z).\]
\end{proposition}

Unfortunately, in this case computing $S_{n}(z)$ is not made any easier through this formula, as computing the functions $S_n^j(z)$ is difficult.  However, this illustrates the general fact that one does not always need to compute the generating function for lattice point counts in every chamber to compute $\bar{H}_G(z)$.  Specifically, when deriving Proposition~\ref{symmetry}, we used the fact that the automorphism group of $K_{1,n-1}$ is $\mathfrak{S}_{n-1}$ and the automorphism action of the group permutes the regions of $\P_n$ by linear transformations.  Thus, we were able to compute the open Ehrhart series for $\P_n$ by calculating the Ehrhart series for one region for each orbit of the group action and sum these series with weights according to the size of each orbit to form the open Ehrhart series for $\P_n$.  When we consider the inside-out polytope $\P_G:=([0,1]^n,\mathcal{H}_0^\sharp (L))$ for the Laplacian $L$ of an arbitrary graph $G$ with automorphism group $\mathrm{Aut}(G)$, the automorphism group will act on the regions and the open Ehrhart polynomial for $\P_G$ will have a decomposition as a weighted sum of Ehrhart polynomials for regions representing orbits of the actions.

Our final goal in this section is to study the vertices of some of the $\P_n^j$'s, yielding insight into the denominator of $S_n(z)$.  First, it is a straightforward exercise to show that $\P_n^1$ is a pyramid over $[0,1]^{n-1}$ and hence all the vertices of $\P_n^1$ are integral. 
Not surprisingly, other regions of $\P_n$ contain rational vertices.  As we are about to see, for the case of $\P_n^2$ the presence of vertex denominators from the set $\{1,2,\ldots,n-1\}$ provides some explanation for the denominators occuring in the rational forms for $S_n(z)$, $n=3,4,5,6$. 

\begin{proposition}For any $n\geq3$, the following are contained in the vertex set of $\P_n^2$.
\[
\left( \begin{array}{c}
\frac{1}{n-1} \\
1\\
0 \\
0 \\
0 \\
\vdots \\
0 \\
0
\end{array} \right), 
\left( \begin{array}{c}
\frac{1}{n-2} \\
1\\
 \frac{1}{n-2} \\
0 \\
0 \\
\vdots \\
0 \\
0
\end{array} \right),
\left( \begin{array}{c}
\frac{1}{n-3} \\
1\\
\frac{1}{n-3} \\
\frac{1}{n-3} \\
0 \\
\vdots \\
0 \\
0
\end{array} \right),
\ldots,
\left( \begin{array}{c}
\frac{1}{2} \\
1\\
\frac{1}{2} \\
\frac{1}{2} \\
\frac{1}{2} \\
\vdots \\
\frac{1}{2} \\
0
\end{array} \right),
\left( \begin{array}{c}
1 \\
1\\
1 \\
1 \\
1 \\
\vdots \\
1 \\
1
\end{array} \right).
\]
\end{proposition}

\begin{proof}
It is straightforward to check that these vectors are contained in $\P_n^2$.  We need only check that these are actually vertices.  The constraint matrices for $\P_n^2$ are of the form
\[\left( \begin{array}{rrrrrrrr}
1-n & 1 & 1 & \cdots & 1 & 1 & \cdots & 1\\
1 & -1 & 0  & \cdots & 0 & 0 & \cdots & 0\\
-1 & 0 & 1 & \cdots & 0 & 0 & \cdots & 0\\
\vdots & \vdots & \vdots & \ddots & \vdots& \vdots & \ddots & \vdots \\
-1 & 0 & 0 & \cdots & 1 & 0 & \cdots & 0 \\
-1  &  0 & 0 & \cdots & 0 & 1 & \cdots & 0 \\
\vdots & \vdots & \vdots & \vdots & \vdots & \vdots & \ddots & \vdots \\
-1 & 0 & 0 & 0 & 0 & 0 & \cdots & 1 \\
\hline
1 & 0 & 0 & 0 & 0 & 0 & \cdots & 0 \\
0 & 1 & 0 & 0 & 0 & 0 & \cdots &0\\
\vdots & \vdots & \vdots & \vdots & \vdots & \vdots & \ddots & \vdots \\
0 & 0 & 0 & 0 & 0 & 0 & \cdots & 1 \\
\hline
-1 & 0 & 0 & 0 & 0 & 0 & \cdots & 0 \\
0 & -1 & 0 & 0 & 0 & 0 & \cdots &0\\
\vdots & \vdots & \vdots & \vdots & \vdots & \vdots & \ddots & \vdots \\
0 & 0 & 0 & 0 & 0 & 0 & \cdots & -1 \\
\end{array} \right)
\left(\begin{array}{c}
x_1\\
x_2\\
\vdots \\
x_n
\end{array}\right)
<
\left( \begin{array}{r}
0 \\
0 \\
0 \\
\vdots \\
0\\
0 \\
\vdots \\
0 \\
\hline
1 \\
1 \\
\vdots \\
1 \\
\hline
0 \\
0 \\
\vdots \\
0 \\
\end{array}\right).\]
For $i=0,\ldots,n-2$, let $v_i$ be the vector $(\frac{1}{n-1-i},1,\frac{1}{n-1-i},\ldots,\frac{1}{n-1-i},0,\ldots,0)^T$.  For $i>0$, the tight inequalities for $v_i$ are given by the first row of the above system, the third through $(3+i-1)\st$ rows of the system, and the rows of the identity matrices where equality holds, i.e., $1$ or $0$ entries in $v_i$.  If $i=0$, the tight inequalities for $v_i$ are given by the first row of the above system and the rows of the identity matrices where equality holds.  The matrix of coefficients for the tight inequalities is thus a square $n\times n$ matrix whose determinant is equal, up to sign, to the $i\times i$ minor $L_{n,i}$ of $L_n$ with rows and columns indexed by the vertices $X_i:=\{1,3,4,\ldots,3+i-1\}$ when $i>0$ and $\{1\}$ when $i=0$.  To show that $v_i$ is a vertex, we need only show that $\det(L_{n,i})\neq 0$.

The non-vanishing of $\det(L_{n,i})$ is trivial when $i=0$.  For the other cases, this follows from standard results in algebraic graph theory, specifically Lemmas 7.4 and 7.5 of \cite{BiggsAGT}.  To summarize the argument, the Cauchy--Binet theorem implies that
\begin{equation}
\det(L_{n,i})=\sum_Y[\det(\partial_1(X_i,Y))]^2, \label{sum}
\end{equation}
where the sum is over all subsets $Y$ of the edges of $K_{1,n-1}$ such that $|X_i|=|Y|$ and $\partial_1(X_i,Y)$ is the submatrix of the matrix corresponding to the boundary map for $K_{1,n-1}$ indexed by the vertices $X_i$ and edges $Y$.  That this sum is non-zero is easily checked by applying \cite[Lemma 7.4]{BiggsAGT}; in the notation of the lemma, to choose a $Y$ that satisfies the conditions of the lemma, we have to include in $Y$ all the edges induced by $X_i$ and we are then allowed to choose one more edge to make $|X_i|=|Y|$.  As we begin with a tree, this choice of $Y$ satisfies the conditions of the lemma.  Hence $\det(\partial_1(X_i,Y))=1$ for some $Y$. As (\ref{sum}) is a sum of squares and one term is nonzero, we have that $\det(L_{n,i})>0$, as desired.
\end{proof}


\bibliographystyle{plain}

\def\cprime{$'$} \def\cprime{$'$}
\providecommand{\bysame}{\leavevmode\hbox to3em{\hrulefill}\thinspace}
\providecommand{\MR}{\relax\ifhmode\unskip\space\fi MR }
\providecommand{\MRhref}[2]{%
  \href{http://www.ams.org/mathscinet-getitem?mr=#1}{#2}
}
\providecommand{\href}[2]{#2}

\bibliography{Braun}

\end{document}

%% file: P3.pstex_t
\begin{picture}(0,0)%
\includegraphics{P3.pstex}%
\end{picture}%
\setlength{\unitlength}{4144sp}%
\begingroup\makeatletter\ifx\SetFigFontNFSS\undefined%
\gdef\SetFigFontNFSS#1#2#3#4#5{%
  \reset@font\fontsize{#1}{#2pt}%
  \fontfamily{#3}\fontseries{#4}\fontshape{#5}%
  \selectfont}%
\fi\endgroup%
\begin{picture}(2299,2082)(284,-886)
\end{picture}%

%% file: Beck_Braun_NonHarmonicColorings_AMS_Final.bbl
\begin{thebibliography}{10}

\bibitem{BeckRobinsCCD}
Matthias Beck and Sinai Robins.
\newblock {\em Computing the continuous discretely}.
\newblock Undergraduate Texts in Mathematics. Springer, New York, 2007.

\bibitem{BeckZasIOP}
Matthias Beck and Thomas Zaslavsky.
\newblock Inside-out polytopes.
\newblock {\em Adv. Math.}, 205(1):134--162, 2006.

\bibitem{BeckZasFlows}
Matthias Beck and Thomas Zaslavsky.
\newblock The number of nowhere-zero flows on graphs and signed graphs.
\newblock {\em J. Combin. Theory Ser. B}, 96(6):901--918, 2006.

\bibitem{BiggsAGT}
Norman Biggs.
\newblock {\em Algebraic graph theory}.
\newblock Cambridge Mathematical Library. Cambridge University Press,
  Cambridge, second edition, 1993.

\bibitem{Birk}
George~D. Birkhoff.
\newblock A determinant formula for the number of ways of coloring a map.
\newblock {\em Ann. of Math. (2)}, 14(1-4):42--46, 1912/13.

\bibitem{breuersanyal}
Felix Breuer and Raman Sanyal.
\newblock Ehrhart theory, modular flow reciprocity, and the {T}utte polynomial.
\newblock http://arxiv.org/abs/0907.0845, to appear in Math. Z.

\bibitem{chenwang}
Beifang Chen and Jue Wang.
\newblock The flow and tension spaces and lattices of signed graphs.
\newblock {\em European J. Combin.}, 30(1):263--279, 2009.

\bibitem{vanherick}
Andrew~Van Herick.
\newblock Theoretical and computational methods for lattice point enumeration
  in inside-out polytopes.
\newblock MA thesis, San Francisco State University, 2007. Available at {\tt
  http://math.sfsu.edu/beck/teach/masters/andrewv.pdf}. \\ {\it Software
  package} {\tt IOP} available at {\tt http://iop.sourceforge.net/}.

\bibitem{kocholflow}
Martin Kochol.
\newblock Polynomials associated with nowhere-zero flows.
\newblock {\em J. Combin. Theory Ser. B}, 84(2):260--269, 2002.

\bibitem{kocholtension}
Martin Kochol.
\newblock Tension polynomials of graphs.
\newblock {\em J. Graph Theory}, 40(3):137--146, 2002.

\bibitem{LaHaye}
Ross La~Haye.
\newblock Binary relations on the power set of an {$n$}-element set.
\newblock {\em J. Integer Seq.}, 12(2):Article 09.2.6, 15, 2009.

\bibitem{LovaszDAFunctions}
L{\'a}szl{\'o} Lov{\'a}sz.
\newblock Discrete analytic functions: an exposition.
\newblock In {\em Surveys in differential geometry. {V}ol. {IX}}, Surv. Differ.
  Geom., IX, pages 241--273. Int. Press, Somerville, MA, 2004.

\bibitem{Soardi}
Paolo~M. Soardi.
\newblock {\em Potential theory on infinite networks}, volume 1590 of {\em
  Lecture Notes in Mathematics}.
\newblock Springer-Verlag, Berlin, 1994.

\bibitem{StanleyAcyclic}
Richard~P. Stanley.
\newblock Acyclic orientations of graphs.
\newblock {\em Discrete Math.}, 5:171--178, 1973.

\bibitem{StanleyVol1}
Richard~P. Stanley.
\newblock {\em Enumerative combinatorics. {V}ol. 1}, volume~49 of {\em
  Cambridge Studies in Advanced Mathematics}.
\newblock Cambridge University Press, Cambridge, 1997.
\newblock With a foreword by Gian-Carlo Rota, Corrected reprint of the 1986
  original.

\bibitem{Whitney}
Hassler Whitney.
\newblock A logical expansion in mathematics.
\newblock {\em Bull. Amer. Math. Soc.}, 38(8):572--579, 1932.

\end{thebibliography}
